\documentclass{amsart}
\usepackage{amssymb,amsmath,amssymb,array,epsfig,graphicx}

\theoremstyle{plain}
\newtheorem{lemma}{Lemma}

\newtheorem{theorem}[lemma]{Theorem}
\newtheorem{corollary}[lemma]{Corollary}

\theoremstyle{definition}
\newtheorem{example}{Example}
\newtheorem{definition}{Definition}

\newcommand{\R}{\mathbb{R}}
\newcommand{\C}{\mathbb{C}}

\newcommand{\spec}{\mathrm{Spec}\,}
\newcommand{\spe}{\mathrm{Spec}}

\newcommand{\be}{\begin{equation}}
\newcommand{\ee}{\end{equation}}
\newcommand{\mat}{\left(\begin{array}}
\newcommand{\rix}{\end{array}\right)}
\newcommand{\eoe}{\null \hfill $\diamond$}

\renewcommand{\qed}{\ifmmode$\Box$\else{\unskip\nobreak\hfil
\penalty50\hskip1em\null\nobreak\hfil$\Box$
\parfillskip=0pt\finalhyphendemerits=0\endgraf}\fi}

\title[]{On the boundary of the pseudospectrum\\
and its fault points}

\author[]{Lyonell Boulton$^{(\ast)}$ and Peter Lancaster}

\date{January 2007}

\thanks{$(\ast)$ Author supported by the Leverhulme Trust.}

\begin{document}
\begin{abstract}
The study of pseudospectra of linear transformations has become a
significant part of numerical linear algebra 
and related areas. A large body of
research activity has focused on how to compute these sets for a
given spectral problem with, possibly, certain underlying
structure. The theme of this paper was motivated by the question: 
How effective are path-following procedures for tracing the pseudospectral 
boundary? The present study of the mathematical properties 
of the boundary of the pseudospectrum is the result. Although this
boundary is generally made up smooth curves, it is shown how the 
Schur triangular form of the matrix can be used to analyse the singular 
points of the boundary.
\end{abstract}

\maketitle

\

\vspace{-.5in}

\section{Preliminaries} \label{s1}

\Large

In this manuscript we discuss regularity properties of the boundary of
the pseudospectrum of a matrix $A\in \C^{n\times n}$. This boundary
turns out to be a piecewise smooth curve. Our main concern
is how the structure of the Schur triangular form of $A$ determines the
singular points on this curve.

Let us begin by adapting the results of \cite{blp} on general matrix
polynomial to the particular case of the linear monic polynomial.
Below, $\spe\,A$ denotes the spectrum of $A$ and $\|\!\cdot\!\|$ denotes
the (maximum) norm of $A$ as a linear operator in the
Euclidean space $\C^n$. For
a given $\delta\geq 0$, the pseudospectrum of $A$ is the set
\[
   \spe_\delta\,A:=\bigcup_{\|E\|\leq \delta} \spe\,(A+E).
\]
By construction $\spe_0\,A=\spe\,A$. For $\delta>0$ sufficiently small,
$\spe_\delta\,A$ consists of ``small'' connected components around
$\spe\,A$. As $\delta$ increases, these components enlarge, collide and
eventually intersect in various complicated ways and the boundary of
$\spe_\delta\,A$, $\partial \spe_\delta\,A$, becomes
more complex. Typically, corners will appear in $\partial \spe_\delta\,A$
as a consequence of two or more of these components intersecting.
As we shall see below, the least singular values of
$P(\lambda)=(\lambda-A)$ plays an important role
in the dynamics of this process.

Let $\lambda=x+iy\equiv(x,y)$ be in the complex plane ($\equiv \R^2$).
The singular values of $P(\lambda)$ are the nonnegative square roots
of the $n$ eigenvalue functions of $P(\lambda)^\ast P(\lambda)$.
They are denoted by
\[
    s_1(\lambda) \geq \ldots \geq s_n(\lambda)\geq 0.
\]
Both the spectrum and the pseudospectra of $A$ are characterised
by the real-valued function $s_n:\C\longrightarrow [0,\infty)$, given
by the smallest singular value. Indeed, it is well known \cite{tre} that 
for all $\delta\geq 0$, 
\begin{equation}\label{sp2}
   \spe_\delta \,A=\{\lambda \in \C\,:\, s_n(\lambda)\leq \delta \}.
\end{equation}
It will be seen that this characterisation is 
crucial for the study of smoothness properties of the boundary, 
$\partial \spe_\delta\,A$. They are the subject of the next
section and the main conclusions are contained in
Theorem \ref{t3}. In particular, the notion of 
{\em fault points} (at which the least singular value is multiple)
is introduced and developed. In Section 3 the Schur triangular form is studied
and it is shown that, by defining a certain equivalence relation on the
eigenvalues of $A$, the general Schur triangle can be reduced to
a {\em block} diagonal form. It is natural to ask whether the
equivalence classes of eigenvalues generated in this way are
independent of the the particular unitary transformation used. This
is shown to be the case when all eigenvalues of $A$ are distinct.

Section 4 includes a classification and study of singular points on the
boundary of a pseudospectrum and, finally, Section 5 is devoted to
a study of these singular points in the case of matrices with size
not exceeding $n=3$.

\section{Regularity properties of the pseudospectral boundary}
\label{ns2}

Let us first consider regularity properties of $s_n(\lambda)$ as a
function defined on the complex plane. Let
\[
   \Sigma_j:=\{(x,y,s_j^2(x,y))\in \R^3 \,:\, (x,y)\in \R^2  \}.
\]
\begin{lemma} \label{t1}
$\bigcup_{j=1}^n\Sigma_j$ is a real algebraic variety.
\end{lemma}
\proof Define the function
\begin{equation}\label{eq.d}
   d(x,y,S)=\mathrm{det}\,\left[ I\,S - P(x+iy)^*P(x+iy)\right],
  \quad (x , y , S) \in \R^3 .
\end{equation}
As the matrix $P(x+iy)^* P(x+iy)$ is hermitian, $d(x,y,S)$
is a polynomial in $\,x,y,S$ with real coefficients, and since
\[
   \bigcup_{j=1}^n\Sigma_j = \{ (x,y,S) \in \R^3 :\,
                                      d(x,y,S) = 0 \} ,
\]
the result follows.\qed

\medskip

Real algebraic surfaces usually have strong smoothness properties, but
self intersections and cusps of various types may also occur as the following
two examples demonstrate.

\smallskip

\begin{example} \label{ex1} 
If $A=\mathrm{diag}\,[1,e^{2i\pi/3},e^{-2i\pi/3}]$, then
\[\spe[P(\lambda)^\ast P(\lambda)]=\{|\lambda-1|^2,|\lambda-e^{2i\pi/3}|^2,
|\lambda-e^{-2i\pi/3}|^2\}.\]
Thus $\Sigma_1\cup \Sigma_2\cup \Sigma_3$ is the union of
three paraboloids of revolution with minima at  $1,e^{2i\pi/3}$ and
$e^{-2i\pi/3}$
respectively. These paraboloids intersect each other in three different planes
parallel to the vertical axis. The only point where the three surfaces
intersect simultaneously is  $(0,0,1)$.   \eoe
\end{example}

\begin{example} \label{ex2} If
\[
   A \,=\, \begin{pmatrix} 3/4 & 1 & 1 \\
                     0 & 5/4 & 1 \\
                     0 & 0 & -3/4   \end{pmatrix} ,
\]
$\Sigma_1\cup\Sigma_2\cup\Sigma_3$ has a conic double point
at $(0,0,5/16)$. At this point $s_3(0)=s_2(0)=\sqrt{5}/4$. \eoe
\end{example}

\smallskip

It will be seen below that the occurrence of isolated
singularities such as that of Example~\ref{ex2} is rare.
In this example, the matrix $A$ had to be carefully crafted to allow
the conic double point around the origin.
Any slight change in the coefficients of $A$ would
eliminate this degeneracy.

It is well known \cite[Theorem~S6.3]{glr} that, for each $j$,
$s_j(\gamma(t))^2$ is a real analytic function of $t\in \R$
whenever $\gamma:\R\longrightarrow \C$ is analytic. Thus,
$s_j(\lambda)=s_k(\lambda)$ for $j\not=k$
on a non-empty open set $\mathcal{O}$ only when
$\mathcal{O}=\C$. Therefore
two different surfaces $\Sigma_j$ can intersect only in a set
of (topological)
dimension at most one (see also \cite{blp}).

To describe these surfaces more precisely,
let $p_A\equiv p:\{1,\ldots,n\}\longrightarrow
\{1,\ldots,n\}$ (not always onto) satisfy the following
properties:
\begin{itemize}
\item[(a)] $\Sigma_j=\Sigma_{p(j)}$.
\item[(b)]
$\bigcup_{j=1}^n \Sigma_j=\bigcup_{j=1}^n \Sigma_{p(j)}$ and
\item[(c)] $\Sigma_{p(j)}=\Sigma_{p(k)} \iff
   p(j)=p(k)$.
\end{itemize}
The map $p$ is a choice of the indexes of the different $\Sigma_j$
and only these.
We denote by $m_A(j)\equiv m(j)=m(p(j))$
the (constant) value of the algebraic multiplicity of
$s_j(\lambda)$ for almost all $\lambda \in \C$.

\smallskip

We are now ready to identify the
region of the plane where
$s_j(\lambda)$ is differentiable. Here and below $\partial_S$
denotes partial differentiation with respect to $S$. With the 
determinant function $d$ as in (\ref{eq.d}), let
\[
    \mathcal{G}_A^j:=\{(x,y)\in \R^2\,:\,
    \partial_S d^{1/m(j)}(x,y,s_j^2(x,y))=0 \}.
\]

\begin{lemma} \label{t2}
Let $\mathcal{O}\subseteq \R^2$ be an open set such that
 $\mathcal{O}\cap
(\mathcal{G}^j_A \cup \spe\,A)=\varnothing$. Then $s_j^2(x,y)$
is differentiable in $\mathcal{O}$.
\end{lemma}
\proof This lemma is a consequence of the implicit function theorem
applied to the function $d(x,y,S)$. Clearly, the polynomial $d(x,y,S)$ is 
continuously differentiable and, since
$\mathcal{O} \cap \spe\,A=\varnothing$, 
$d(x,y,s_j^2(x,y))\not=0$ for $(x,y)\in \mathcal{O}$.
These two facts ensure that $d^{1/m(j)}(x,y,S)$ is also continuously
differentiable in a suitable $\R^3$ neighbourhood of $\Sigma_j\cap
(\mathcal{O} \times \R)$. Since $\mathcal{O} \cap \mathcal{G}^j_A=\varnothing$,
the hypotheses of the implicit function theorem
applied to $\Sigma_j$ are fulfilled. Since $s_j^2(x,y)$ is the actual
``implicit'' function given locally (and hence globally) by
this theorem, the result follows.   \qed

Note that we can compute $\nabla s_j^2$ on $\mathcal{O}$
by implicit differentiation.

The sets $\mathcal{G}^j_A$ introduced above are closely related to the
set of ``fault points'' introduced in \cite{blp}. For $j=1,\ldots,n$, let
\[
    \mathcal{F}^j_A:=\{(x,y):\, s_{p(j)}(x,y)=s_{p(k)}(x,y)\
    \mathrm{for\ some}\ k\ \mathrm{s.t.}\ p(j)\not=p(k) \}.
\]
Here, we call $\mathcal{F}_A:=\mathcal{F}^n_A$ the set of {\em fault 
points of} $A$. Typically  $\mathcal{F}_A$ is made up of fault {\em lines
or curves} where $\nabla s_n(\lambda)$ is undefined.

Let $\{p_1,\ldots,p_l\}:=p(\{1,\ldots,n\})$. By construction,
$s_{p_j}\not = s_{p_k}$ almost everywhere and
\[
    d(x,y,S)=(S-s_{p_1}^2)^{m(p_1)}\cdots (S-s_{p_l}^2)^{m(p_l)}.
\]
Also
\[
   \partial_S d^{1/m(p_i)}(x,y,S)=
    \sum_{j=1}^l \frac{m(p_j)}{m(p_i)}
(S-s_{p_j}^2)^{\frac{m(p_j)}{m(p_i)}-1}
\prod _{k\not=j,\,  1\leq k\leq n} (S-s^2_{p_k})^{\frac{m(p_k)}{m(p_i)}}.
\]
So $(\tilde{x},\tilde{y})$ satisfies $s_{p_i}(\tilde{x},\tilde{y})
\not= s_{p_j}(\tilde{x},\tilde{y})$ for all $i\not=j$ if and only if
\[
   \partial_S d^{1/m(p_i)}(\tilde{x},\tilde{y},s_{p_i}(\tilde{x},\tilde{y}))
\not = 0.
\]
It  follows from this observation that
\begin{equation} \label{e2}
    \mathcal{G}^j_A = \mathcal{F}^j_A, \qquad
j=1,\ldots,n.
\end{equation}
In particular $s_n(\lambda)$ will be differentiable outside
$\mathcal{F_A}$. This justifies the name chosen for the latter set.

\medskip

Clearly \cite{blp},
\begin{equation} \label{e3}
    \partial \spe_\delta\, A\subseteq \{z\in \C\,:\,
      s_n(\lambda)=\delta\}.
\end{equation}
Note that equality in \eqref{e3} does
not hold in general.
In Example 1 above the origin is in the right hand set, but it is an
interior point of $\spe_1\,A$.

\begin{theorem} \label{t3}
For any $\delta>0$ the boundary of
$\spe\,_{\delta}\,A$ is a piecewise smooth portion of
an algebraic curve. In particular, it has a finite number of
singularities. These singularities are either cusps or
self-intersections. If
$\lambda_0\in\partial \spe\,_{\delta}\,A$ is
a cusp, then $\lambda_0\in \mathcal{F}_A$.
\end{theorem}

\proof The first part of the theorem is Theorem 7 of \cite{blp}.
For the latter part note that if a singularity occurs at
$\lambda_0\in \partial \spe\,_{\delta}\,A$, then
either $\lambda_0\in \mathcal{F}_A$ or $\nabla s_n(\lambda_0)=0$.
If $\nabla s_n(\lambda_0)=0$, then either $\lambda_0$ is a point of 
self-intersection or $\lambda_0\in \spe\, A$. Since no eigenvalue of $A$
lies on the pseudospectral boundary, the result follows.\qed

Theorem~\ref{t3} extends to matrix polynomials without much effort.

\section{Refinement of the Schur triangular form} 
\label{s2}

We denote a Schur factorisation of $A$ by $A=USU^\ast$,  where
\be \label{eq.t}
   S=\begin{pmatrix}  \alpha_{1}&t_{12}&\ldots & t_{1n} \\
   0 & \alpha_2& \ddots   & \vdots \\ \vdots
    & \ddots & \ddots & t_{(n-1)n} \\  0 & \ldots & 0 & \alpha_n
    \end{pmatrix}
\ee
and $U$ is a unitary matrix.
Note that $\spe \,A$ consists of the distinct $\alpha_j$, and that
$S$ and $U$ are not uniquely defined. The singular values
of $P(\lambda)$ are invariant under unitary similarity transformation,
so the Schur form
is an invariant as far as spectrum, pseudospectrum and fault points
are concerned.

The characterisation of all possible Schur triangular forms
of a general matrix $A$ is certainly beyond the scope of this manuscript.
However, Theorems~\ref{t8} and~\ref{t7} below suggest how this issue can
be approached. Some preliminary considerations are required.

\begin{definition} \label{t4}
Let $A=U^\ast SU$ be  a Schur factorisation of $A$ with $S$ as in 
(\ref{eq.t}):\\
(a) We write $\alpha_i \bumpeq \alpha_j$ and  say that
$\alpha_i$ is \emph{directly related} to $\alpha_j$
if one of the following conditions holds:
\[  \textup{(i)} \;\, i=j,\quad
     \textup{(ii)} \;\, i<j \quad \textup{and} \quad t_{ij}\ne 0, \quad
 \textup{(iii)} \;\, i>j \quad \textup{and} \quad t_{ji}\ne 0. \]
(b) We write $\alpha_i\Bumpeq \alpha_j$ and say
that $\alpha_i$ is \emph{block equivalent} 
to $\alpha_j$ if and only if there exists
a subset $\{\sigma_1,\ldots,\sigma_m\}$ of $\{1,\ldots,n \}$ such that
\[
\alpha_i \bumpeq \alpha_{\sigma_1}\bumpeq
\cdots \bumpeq \alpha_{\sigma_m}\bumpeq \alpha_{j}.
\]
\end{definition}

In this definition we consider that any two eigenvalues
on the diagonal of $S$ are ``different'' even in the presence of multiplicity.
The symbol ``$\Bumpeq$'' obviously defines an equivalence relation on 
the set $\{\alpha_i\}_{i=1}^n$.

Let $D$ be the binary matrix obtained from $S$ by preserving the
zeros of the latter and replacing the non-zero entries by $1$. Then $D$
is the adjacency matrix (see \cite{fiedler}) of a graph, $G$.
Two diagonal entries of $S$ will be block equivalent,
$\alpha_i \Bumpeq \alpha_j$, if and only if 
the nodes $i$ and $j$ of $G$ are connected with
a path.  

\begin{lemma} \label{t5}
Let $\sigma$ be any permutation of $\{1,\ldots,n\}$ and $S$ be
upper triangular, as in \textup{(\ref{eq.t})}. Then there is a
unitary matrix  $V$ such that
\[
   V^\ast S V=\begin{pmatrix}  \alpha_{\sigma(1)}&r_{12}&\ldots & r_{1n} \\
   0 & \alpha_{\sigma(2)}& \ddots   & \vdots \\ \vdots
    & \ddots & \ddots & r_{(n-1)n} \\  0 & \ldots & 0 & \alpha_{\sigma(n)}
    \end{pmatrix}
\]
and $\alpha_i\Bumpeq \alpha_j$ if and only if 
$\alpha_{\sigma(i)}\Bumpeq \alpha_{\sigma(j)}$.
\end{lemma}
\proof  By writing the permutation $\sigma$ as a product of
transpositions (cycles of length 2), we see that the case $n=2$ can play
an important role. In this case
$S=\begin{pmatrix} \alpha_1 & t \\ 0& \alpha_2 \end{pmatrix}$ 
with $t\ne 0$ (so that $\alpha_1 \Bumpeq \alpha_2$) and there is just one 
permutation $\sigma$ of interest: $(1,2)\rightarrow (2,1)$. 

Let $\alpha_2-\alpha_1 =be^{i\theta}$ $(b\ge 0)$, and $a=|t|$ and
consider the real orthogonal matrix
\[  W:=\frac{1}{\sqrt{a^2+b^2}}\mat{cc} a & -b\\
                                       b & a\rix. \]
A little calculation shows that 
\[     W^*SW=\mat{cc} \alpha_2 & |t|e^{i\theta}\\
                          0 & \alpha_1 \rix.  \]
Since $t\ne 0$, $|t|e^{i\theta}\ne 0$, and  $\alpha_2 \Bumpeq \alpha_1$,
as required.

For the case $n>2$ let $a_k=|t_{k,k+1}|$, $\alpha_{k+1}-\alpha_k =
 b_ke^{i\theta_k},$ ($b_k\ge 0$), and let
\[    W_k=\frac{1}{\sqrt{a_k^2+b_k^2}}\mat{cc} a_k & -b_k\\
                       b_k & a_k \rix,   \]
Define
\[
   V_k=\mathrm{diag}\;(I_{k-1},\,W_k,\;I_{n-k-1}).
\]
Using a conforming block structure for $S$, 
\[
   V_k^\ast SV_k=\left( \begin{array}{ccc} S_{11}& S_{12}W_k & S_{13}  \\ 
 0 & W_k^*S_{22}W_k & W_k^*S_{23}  \\ 
 0 & 0 & S_{33} \end{array} \right)=:R
\]
and 
\[
  W_k^*S_{22}W_k =\begin{pmatrix} \alpha_{k+1} & 
t_{k+1,k}e^{i\theta_k} \\ 0 & \alpha_k \end{pmatrix}.
\] 
By construction, the upper-right entry of the latter matrix is zero
if and only if $t_{k+1,k}=0$.
Thus $\alpha_k\bumpeq \alpha_{k+1}$ (using $S$) if and only if
$\alpha_{k+1}\bumpeq \alpha_{k}$ (using $R$). 

Since the general permutation $\sigma$ can be expressed as a product
of ``contiguous'' cycles as used above, the required transforming matrix
$V$ can be expressed as a product of the elementary unitary matrices
$V_k$ introduced above.
\qed

\medskip

We now show that a block diagonal factorisation
for $A$ can be obtained from the relation ``$\Bumpeq$''. 

\begin{theorem} \label{t8}
Given any upper triangular matrix $S \in \C^{n\times n}$ there exists 
a unitary matrix $U$ preserving ``$\Bumpeq$'' such that 
\begin{equation} \label{e4}
    U^\ast SU=\mathrm{diag}\,[B_1,\ldots,B_k],
\end{equation}
where, for $l=1,\ldots,k$, 
\[
    B_l=\begin{pmatrix}  \tilde{\alpha}_{1}(l)&r_{12}(l)&\ldots & r_{1m_l}(l) \\
   0 & \tilde{\alpha}_2(l)& \ddots   & \vdots \\ \vdots
&\ddots & \ddots & r_{(m_l-1)m_l}(l) \\ 0 & \ldots & 0 &\tilde{\alpha}_{m_l}(l)
    \end{pmatrix}, 
\]
and
$\tilde{\alpha}_i(l)\Bumpeq \tilde{\alpha}_j(l)$ for all $i,j$ and $l$.
\end{theorem}

\proof
Since ``$\Bumpeq$'' is an equivalence relation, it partitions $\spec S$ into
a family of equivalence classes. So there is a permutation

\[ \sigma =\{\sigma(1),\ldots ,\sigma(n)\} \] 
and there are positive integers
$l_1<l_2<\ldots <l_k<n$ such that
\[\{
  \{\alpha_{\sigma(1)},\ldots,\alpha_{\sigma(l_1)}\},
  \{\alpha_{\sigma(l_1+1)},\ldots,\alpha_{\sigma(l_2)}\},
  \ldots,
  \{\alpha_{\sigma(l_k +1)},\ldots,\alpha_{\sigma(n)}\}\}
\]
is the partition of $\spec S$ under ``$\Bumpeq$''.
The $U$ found in Lemma~\ref{t5}
for this $\sigma$ provides the desired conclusion. \qed

\medskip

In this theorem, $\tilde{\alpha}_i(l) \Bumpeq \tilde{\alpha}_j(\tilde{l})$
if and only if $l = \tilde{l}$. The blocks $B_l$ are not
necessarily unique  in the representation
$U^\ast SU=\mathrm{diag}\,[B_1,\ldots,B_k]$. The
diagonal blocks can be permuted, for example.

Note that the equivalence classes given by ``$\Bumpeq$'' 
on the set $\{\alpha_i\}_{i=1}^n$, correspond to
the connected components of the graph $G$ introduced above. 
When $S$ is sparse, strategies for computing the 
block diagonalisation \eqref{e4} may be based on this observation.

\medskip

It is natural to ask whether Definition~\ref{t4} is independent 
of $U$ and $S$ in the Schur decomposition for $A$. We obtain a positive 
answer if all eigenvalues of $A$ are simple.

\begin{theorem}  \label{t7}
Suppose that $A$ has no multiple eigenvalues, i.e. $\alpha_i\not = \alpha_j$ 
for $i\not=j$. Then the equivalence classes determined by the
equivalence relation ``$\Bumpeq$'' of Definition~\ref{t4} are
independent of the Schur factorisation chosen for $A$.
\end{theorem}
\proof This result is a direct consequence of the following observation.
If $S$ and $T$ are two triangular matrices with the same diagonal entries
$s_{ii}=t_{ii}=\alpha_i$ and with $\alpha_i\not=\alpha_j$ for $i\not=j$,
such that $U^\ast TU=S$ for a suitable unitary matrix $U$,
then $U$ must be diagonal. See \cite[Theorem~2.3]{sha} or \cite{lit}.
\qed
\medskip

The results reviewed at length in \cite{sha} suggest that extension
of this theorem to admit multiple eigenvalues of $A$ would be difficult.

\section{The singular points on the pseudospectral boundary}
\label{ns4}

As we shall see next, the block diagonalisation found 
in Theorem~\ref{t8} provides a natural
classification of the singular points on the boundary of the
pseudospectra of $A$. 

Let $S$ be a Schur triangular form of $A$ and
$B_l$ be the blocks corresponding to 
the diagonalisation of $S$ given in $\eqref{e4}$. 
It is clear from the definition that the singular values of 
$P(\lambda)=\lambda - A $ are those of $(\lambda-B_l)$ for 
$l=1,\ldots,k$. Then it follows from (\ref{sp2}) that
\begin{equation} \label{e1}
    \spe_\delta A=\bigcup _{l=1}^k \spe_\delta B_l,
    \qquad \mathrm{for\ all}\ \delta\geq 0.
\end{equation}

This decomposition motivates the following classification of
singular points on the boundary of the pseudospectrum:

\begin{definition} \label{d2}
We say that $\lambda_0 \in \partial \spe_{\delta}\,A$ is a
\begin{itemize}
\item[-] \emph{stationary point} if $\lambda_0\not\in \mathcal{F}_A$ but
$\nabla s_n(\lambda_0)=0$,
\item[-] \emph{essential fault point} if $\lambda_0\in\mathcal{F}_{B_l}$ 
for some $l=1,\ldots,k$,
 \item[-] \emph{regular fault point} if $\lambda_0\in
\mathcal{F}_A$ but it is not an essential fault points.
\end{itemize}
\end{definition}

Clearly a stationary point of $A$ must also be a stationary point
of $B_l$ for some (not necessary unique) $l=1,\ldots,k$.

If $\lambda_0\in \mathcal{F}_A$, then either
$\lambda_0\in \bigcup_{l=1}^k \mathcal{F}_{B_l}$ or
\[
 \lambda_0\in (\partial\spe_{\delta}B_k\, \cap\,
 \partial\spe_{\delta}B_l)\setminus \bigcup_{l=1}^k \mathcal{F}_{B_l}
\]
for some $k\not=l$. In the former case, $\lambda_0$ is an essential 
fault point and in the latter it is a regular fault point.

Since it is formed as a consequence of two pseudospectra of different blocks 
$B_l$ intersecting, whenever non-empty,
the regular portion of $\mathcal{F}_A$ is expected to be of topological 
dimension $1$. On the other hand, we shall argue in Section~\ref{ns5}
that blocks of small size have only a limited number of essential
fault points. 

\medskip

The role played by the classification just introduced in the dynamics of the
pseudospectrum of $A$ can be better visualised by means of concrete examples.
However, let us first establish two elementary consequences of \eqref{e1}. 

\begin{corollary} \label{t12}
If $A$ is a normal matrix, then all singular points on 
$\partial \spe _\delta A$ are regular fault points. Furthermore,
$\mathcal{F}_A$ is the Voronoi diagram associated to $\spe\,A$.
\end{corollary}
\proof The proof is straightforward. 
See Example~\ref{ex1} for illustration. \qed

\medskip

\begin{corollary} \label{t13}
Let $A$ be unitarily similar to a bi-diagonal matrix. 
Then $\mathcal{F}_A=\varnothing$ if and 
only if no entry in the off-diagonal of $A$ vanishes.
\end{corollary}
\proof Use the fact that $(\lambda-B)^\ast(\lambda-B)$ 
is a tri-diagonal symmetric
matrix. See \cite[\S5.36]{wil}. \qed

\medskip

Corollary {\ref{t13} implies that bi-diagonal matrices 
have no essential fault points.

\begin{example} \label{ex3} A detailed analysis is made of a problem with
two diagonal blocks ($k=2$); see also Figure~\ref{f1}. Let
\[
   A=\begin{pmatrix} 3 & 0 & 0 \\ 0 & -1 & 1 \\ 0 & 0 & 1 \end{pmatrix}.
\]
Any pseudospectrum of $A$ is the union of those of the blocks
\[
     B_1=3 \qquad \mathrm{and} \qquad B_2=\begin{pmatrix}
     -1 & 1 \\ 0 & 1 \end{pmatrix}.
\]
The (only) singular value of $(x+iy-B_1)$ is $\sqrt{(x-3)^2+y^2}$,
so $\spe_\delta B_1$ is a disk centred at $(x,y)=(3,0)$ of radius
$\delta$.
The singular values of $B_2$ are
\[
s_{\pm}(x+iy)=\left(\frac{3}{2}+x^2+y^2\pm \frac{\sqrt{5+20x^2+4y^2}}{2}\right)
^{1/2}.
\]
It is straightforward to see that the least singular value
$s_-(x+iy)=\delta$ if and only if
\[
   (x^2+y^2)^2-2(1+\delta^2)x^2+2(1-\delta^2)y^2+1-3\delta^2
+\delta^4 =0.
\]
This shows that $\spe_\delta B_2$ is a spiric section 
for all
$\delta>0$. It is also straightforward to see that
$\mathcal{F}_A$ is  a portion of the hyperbola given explicitly by
\[
    \mathcal{F}_A=\{(x,y)\,:\, 31\,x^2-y^2-90\,x+55=0, x>1 \}.
\]
Any singularity on $\partial \spe_\delta\,A$
which arises as a consequence of
$\partial \spe_\delta\,B_1$ intersecting with
$\partial \spe_\delta\,B_2$, will be  a regular fault point.
Also, the curve
$\partial \spe_{\delta}\,B_2$  has a stationary point at the
origin when $\delta=\sqrt{\frac{3-\sqrt{5}}{2}}$. 

Note that
a further self-intersection occurs on $\partial \spe_\delta\,A$
at $(x,y)=(\frac{45+8\sqrt{5}}{31},0)$ for
$\delta=\frac{48-8\sqrt{5}}{31}$ since the boundaries of 
$\spe_\delta\,B_1$ and $\spe_\delta\,B_2$
touch at this point. The pseudospectrum of $A$ will consist of
three connected components for $0\leq \delta<\sqrt{\frac{3-\sqrt{5}}{2}}$,
two components for $\sqrt{\frac{3-\sqrt{5}}{2}}\leq
\delta<\frac{48-8\sqrt{5}}{31}$ and a single component
for  all  $\delta\geq\frac{48-8\sqrt{5}}{31}$. 
\eoe  \end{example}

\medskip

The next observation may be relevant in
the effective design of a corrector step in path-following
algorithms for tracking $\partial \spe_\delta\, A$. 
The second part asserts that the corners in 
$\partial \spe_\delta\, A$ at regular fault points will always
be re-entrant. 

\begin{lemma} \label{t11}
Let $\lambda_0$ be a regular fault point in $\partial \spe_\delta\, A$.
Assume that the curve  $\partial \spe_\delta\, A$ fails to have
a tangent line at $\lambda_0$ and does not self-intersect
at this point.
Then we can always find an $\alpha > \pi$, depending only on $\lambda_0$, 
satisfying the following property: If $0<\beta < \alpha$, there exists
$r>0$ and $0\leq \gamma<2\pi$ such that the sector
\[
   \{\lambda=\lambda_0+\rho e^{i(\theta+\gamma)}\,:\, 0\leq \theta \leq \beta,
\,0\leq \rho \leq r \} \subset \spe_\delta\,A. 
\]
\end{lemma}
\proof This is a consequence of the fact that a 
regular fault point of $\partial \spe_\delta\, A$ is formed 
by two different intersecting connected
components. \qed

\medskip 

 \begin{theorem}  \label{t9}
If $\alpha_i,\, \alpha_j \in \;\textup{Spec}A$ are in the same 
connected component of $\C \setminus \mathcal{F}_{A}$, then 
$\alpha_i\Bumpeq \alpha_j$.
\end{theorem}
\proof
It suffices to show that if $\alpha_i\not\Bumpeq \alpha_j$, then
any continuous trajectory $\phi:[0,1]\longrightarrow \C$ such
that $\phi(0)=\alpha_i$ and $\phi(1)=\alpha_j$, intersects $\mathcal{F}_{S}$
(see \eqref{eq.t}).
We achieve this by noticing that $\alpha_i\in \mathrm{diag}[B_l]$
and $\alpha_j \in \mathrm{diag} [B_m]$ for $l\not=m$, and applying
the mean value theorem inductively to $\rho_i-\rho_j$, where
\[
   \rho_k(t):=\mathrm{min.sing.val}(\phi(t)-B_p), \qquad \alpha_k\in
\mathrm{diag}\,[B_p].
\]
\qed

The converse of Theorem~\ref{t9} does not
hold in general. It is easy to construct  examples where
$\alpha_i\Bumpeq \alpha_j$, but $\alpha_i$, $\alpha_j$
both belong to different components of $\C\setminus \mathcal{F}_{A}$.
One such example is the following.

\begin{example} \label{ex5} If
\[
  A=\left(\begin{array}{cccc}-1 & 1 & 0 & 0 \\0 & 1 & 0 & 0
\\0 & 0 & -i & 1 \\0 & 0 & 0 & i\end{array}\right),
\]
$\mathcal{F}_{A}=\{\lambda \in \C : \mathrm{Re}\,
\lambda=\pm \mathrm{Im}\,\lambda\}$.
\eoe  \end{example}

\medskip

\section{Essential fault points and the Schur structure of  small matrices}
\label{ns5}

Determining the structure of $\mathcal{F}_A$ for a given matrix $A$,
is typically involved. As we confirmed in the previous section,
the set of fault points
can be empty or consist of a single point, but
it can also be a complicated set such as a
Voronoi diagram. By virtue of \eqref{e1}, the set of regular fault points
is completely characterised  once the pseudospectra of each of the blocks
in the diagonal factorisation \eqref{e4} are known. 
In this section we show that, for sufficiently small matrices, the number of
essential fault points is finite.

It is easy to characterise 
the set of fault points of a $2\times 2$ 
triangular matrix
\[
       B=\begin{pmatrix} \alpha & r \\ 0 & \beta \end{pmatrix}.
\]
If $\alpha \Bumpeq \beta$, then $\mathcal{F}_B=\varnothing$. 
If $\alpha\not\Bumpeq \beta$, 
then 
\[
   \mathcal{F}_B=\{\lambda\in\C\,:\, |\lambda-\alpha|=|\lambda-\beta| \}
\]
for $\alpha\not=\beta$ and
$\mathcal{F}_B=\varnothing$ for $\alpha=\beta$. 
By virtue of Corollary~\ref{t13}, no fault point of a $2\times 2$ matrix 
can be essential.

\medskip

The $3\times 3$ case is more involved. This is illustrated in the following 
example (which is a generalisation of  Example~\ref{ex2}).

\begin{example} \label{ex4} See Figure~\ref{f2}. Let
\[
   A=\begin{pmatrix} a &1 &1 \\ 0 & 5/4 & 1 \\ 0 & 0 & c \end{pmatrix},
\]
where $|a|=|c|=3/4$. Then $\mathcal{F}_A=\{0\}$ and $a\Bumpeq 5/4 \Bumpeq c$. Whenever $0\in
\partial \spe_\delta\,A$, it will be an essential fault point.
This can only occur at $\delta=s_3(0)=\sqrt{5}/4$.

The complicated dynamic of
the pseudospectral boundary as we move the
parameters $a$ and $c$ on the circle with radius $3/4$, is
illustrated in  Figure~\ref{f2}. There we show the evolution of
the essential singularity when $a=\overline{c}=3e^{i\theta}/4$
for $\theta=k\pi/4$, $k=0,1,2,3,4$. When $k=0$, 
$\spe_{\sqrt{5}/4}A$ is connected and there is clear indication
of a re-entrant corner on the left side of the boundary. 
The re-entrant angle at this corner becomes more
pronounced when $k=1$.
For $k=2,3$ the concavity prevails. However there exist 
critical $\theta_1\in(\pi/4,\pi/2)$ and $\theta_2\in(3\pi/4,\pi)$, such that
$0\not \in \partial\spe_{\sqrt{5}/4}A$ for $\theta\in(\theta_1,\theta_2)$,
and hence there is no corner on $\partial\spe_{\sqrt{5}/4}A$. 
At the final stage $k=4$, the pseudospectrum now consists of two connected components,
the concavity on the left side has vanished, however a new
re-entrant corner forms in a different part of the boundary.
\eoe  \end{example}

\medskip

More generally, let
\begin{equation} \label{e5}
       B=\begin{pmatrix} \alpha_1 & r & s \\ 0 & \alpha_2 & t
       \\ 0&0&\alpha_3 \end{pmatrix}.
\end{equation}
When $\alpha_j \not \Bumpeq \alpha_k$ for some $j\not=k$,
 $B$ can be reduced to 
a block diagonal form of smaller size, so the set of fault points
is characterised by using \eqref{e1}.
If $\alpha_1 \Bumpeq \alpha_2 \Bumpeq \alpha_3$,
every point in $\mathcal{F}_B$ is an essential fault point and we have 
the following result.

\begin{theorem} \label{t14}
Let $B$ be as in \eqref{e5} and assume that no two of the upper 
triangular entries $r,s,t$ 
vanish simultaneously (i.e. $\alpha_1\Bumpeq \alpha_2\Bumpeq
\alpha_3$). If $rst=0$, then $\mathcal{F}_B=\varnothing$.
If $rst\not=0$, then $\mathcal{F}_B$ is either empty or 
consists of a single point. In the latter case,
\begin{equation} \label{ee5}
    \mathcal{F}_B\subset \{\lambda\in\C\,:\, 
\arg (\alpha_2-\lambda)=\arg (r\overline{s}t)\}.
\end{equation}
\end{theorem}

\proof If $rst=0$, the claimed assertion is a consequence of 
Corollary~\ref{t13} and a suitable permutation of the rows and columns of
$B$.

Let $a=(\alpha_1-\lambda)$, $b=(\alpha_2-\lambda)$, $c=(\alpha_3-\lambda)$,
and assume that $rst\not=0$. Let 
\[
B_1=\begin{pmatrix} a & r \\ 0 & b \end{pmatrix} \qquad 
\mathrm{and}
\qquad
u=\begin{pmatrix} s \\ t \end{pmatrix},
\]
so that
\[
   B-\lambda=\left(\begin{array}{c|c} B_1 & u \\ \hline 0 & c \end{array} 
   \right),
\]
and
\[
   (B-\lambda)^\ast(B-\lambda)=\left(\begin{array}{c|c} B_1^\ast B_1 & v \\ 
\hline v^\ast & \|u^2\|+|c|^2 \end{array} 
   \right),
\]
where
$v=B_1^\ast u$. By the Cauchy interlacing theorem, if $B^\ast B$ has a double 
eigenvalue $\sigma=s_3^2=s_2^2$, then $\sigma$ is also the 
minimal eigenvalue of $B_1^\ast B_1$. 

Since $r\not=0$,
$B_1^\ast B_1$ has only simple eigenvalues. Let $0\not=e\in \C^2$ be an 
eigenvector such that  $B_1^\ast B_1 e= \sigma e$. A straightforward argument
shows that, if $\sigma$ is a double eigenvalue of $B^\ast B$, then
$v\perp e$.

Now,  $v=\begin{pmatrix} \overline{a}s  \\ \overline{r} s+\overline{b}t 
   \end{pmatrix}$, 
\begin{gather*}
2\sigma=|a|^2+|b|^2+|r|^2-\sqrt{(|a|^2+|b|^2+|r|^2)^2-4|a|^2|b|^2} \qquad
\mathrm{and} \\ e=\begin{pmatrix} 
|b|^2-|a|^2+|r|^2-\sqrt{(|a|^2+|b|^2+|r|^2)^2-4|a|^2|b|^2} \\
-2a\overline{r} \end{pmatrix}.
\end{gather*}
Note that all these three quantities depend on $\lambda$. If $v^\ast e=0$, 
then
\[
    q(\lambda)=\overline{s}(|b|^2-|a|^2-|r|^2-
\sqrt{(|a|^2+|b|^2+|r|^2)^2-4|a|^2|b|^2})-2\overline{r}\overline{t}b=0.
\]
As the coefficient of $\overline{s}$ in the above expression is real,
\eqref{ee5} is guaranteed. 

Let $\gamma(t)=\alpha_2+te^{i\arg (r\overline{s}t)}$, 
$-\infty<t<\infty$, be a parameterisation of the line where $\mathcal{F}_A$
lies. Then $q(\gamma(t))=q_1(t)-\sqrt{q_2(t)}$ where $q_1$ depends linearly
in $t$ and $q_2$ is a quadratic polynomial in $t$. Moreover, $(q_1(t))^2$ and
$q_2(t)$ have the same coefficient of order $2$ in $t$. 
Thus $q(\gamma(t))$ can 
only vanish at no more than one $t$-value, $t_1$. Since 
$\gamma(t_1)$ is the only possible 
essential fault point of $B$, the proof is complete.  \qed  

\medskip

In Example~\ref{ex4} the line described by 
the right hand expression in \eqref{ee5} is the real axis and $\mathcal{F}_A$
is the origin.

It is natural to expect that the argument presented in
the proof of Theorem~\ref{t14} can be extended inductively
to blocks of larger size. We have not explored this possibility in much
detail. However our observations lead us to conjecture that, perhaps,
the number of essential fault points is always finite for matrices of any 
size. This issue certainly requires further investigation.



\newpage

\begin{figure}[t]
\epsfig{file=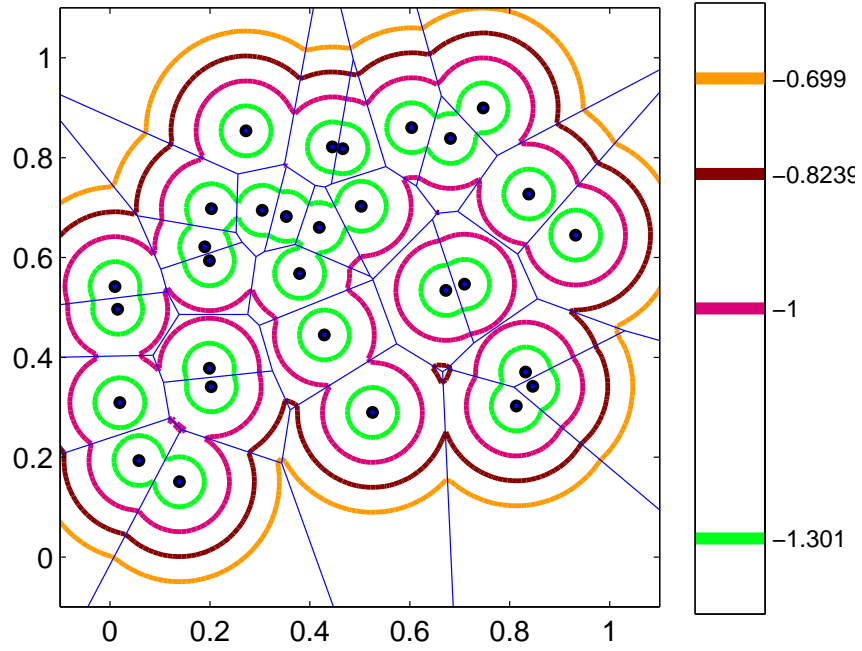, height=10cm} 
\caption{Pseudospectra of a randomly generated diagonal matrix
along with the Voronoi diagram associated to the spectrum. \label{f3}}
\end{figure}

\begin{figure}[t]
\epsfig{file=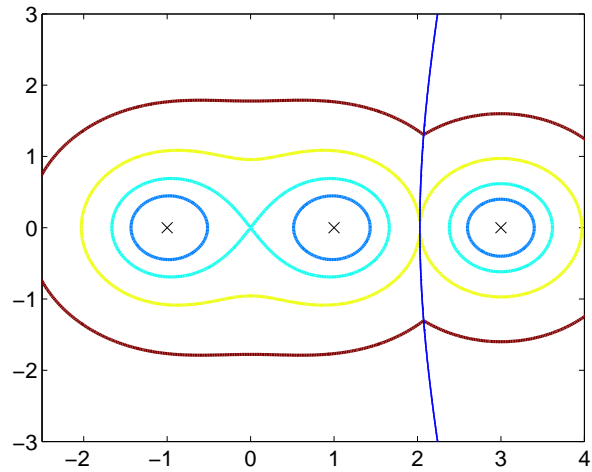, height=7cm} \caption{Example~\ref{ex3}.
The thick lines are $\partial \spe_\delta\,A$ for $\delta=2/5,
\sqrt{\frac{3-\sqrt{5}}{2}},\frac{48-8\sqrt{5}}{31},8/5$. The thin line
is $\mathcal{F}_A$. \label{f1}}
\end{figure}

\newpage

\begin{figure}
\epsfig{file=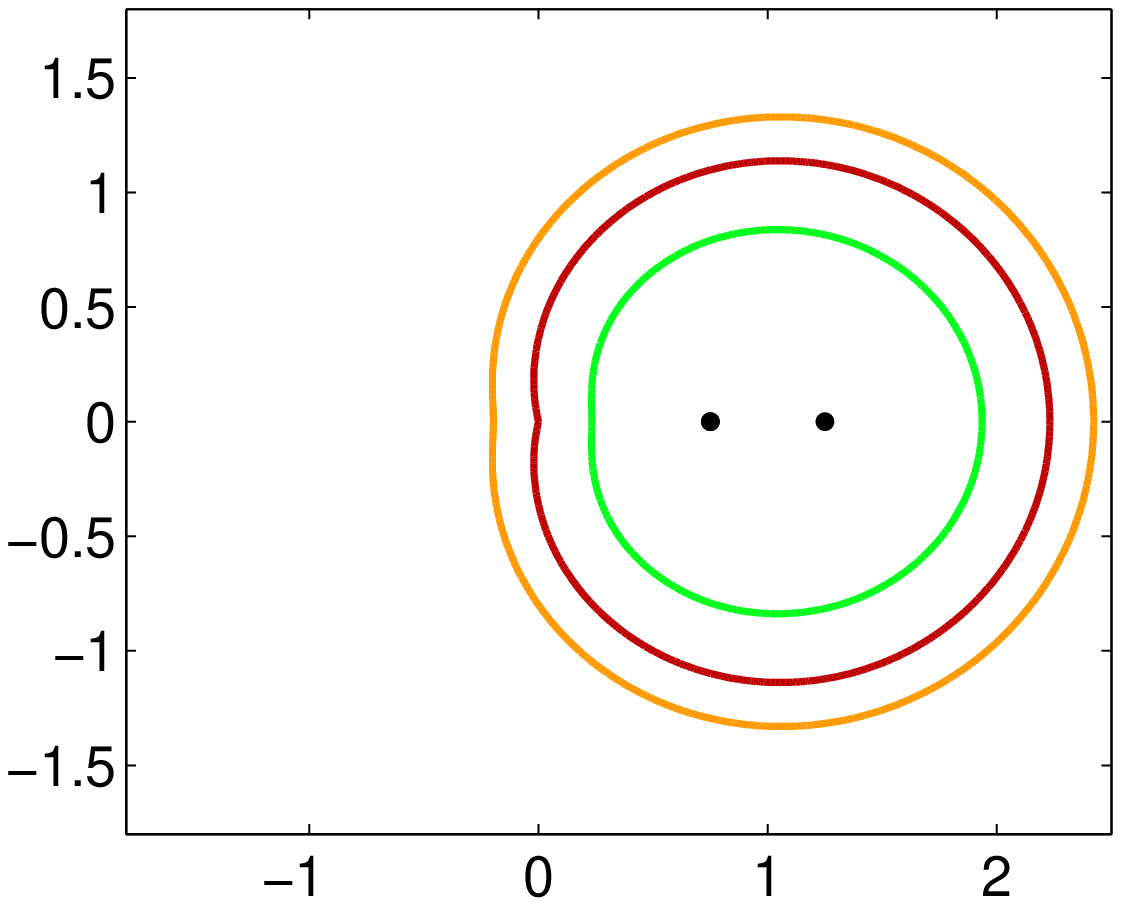, height=4cm}\hspace{-2cm}
\epsfig{file=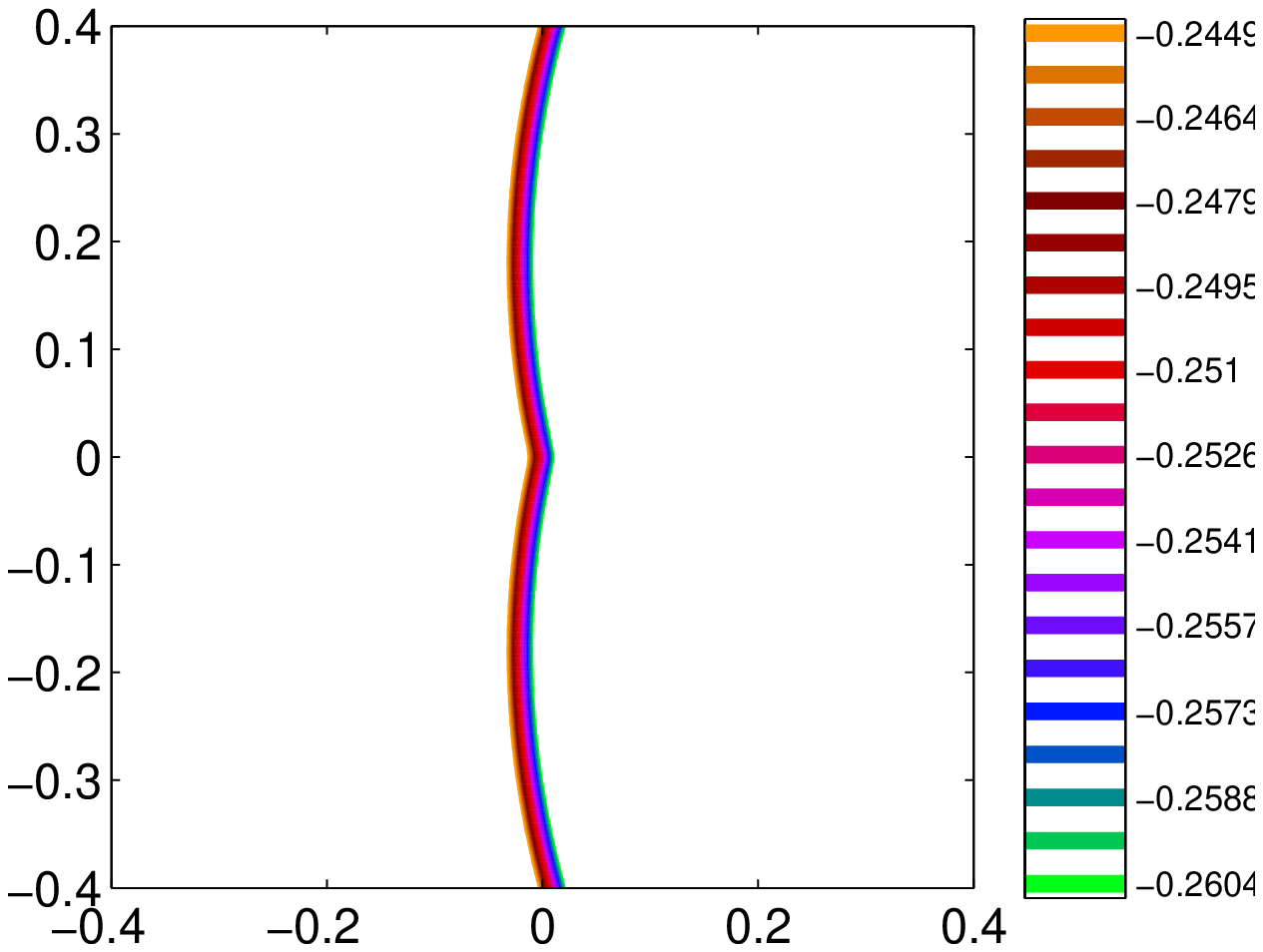, height=4cm}

\epsfig{file=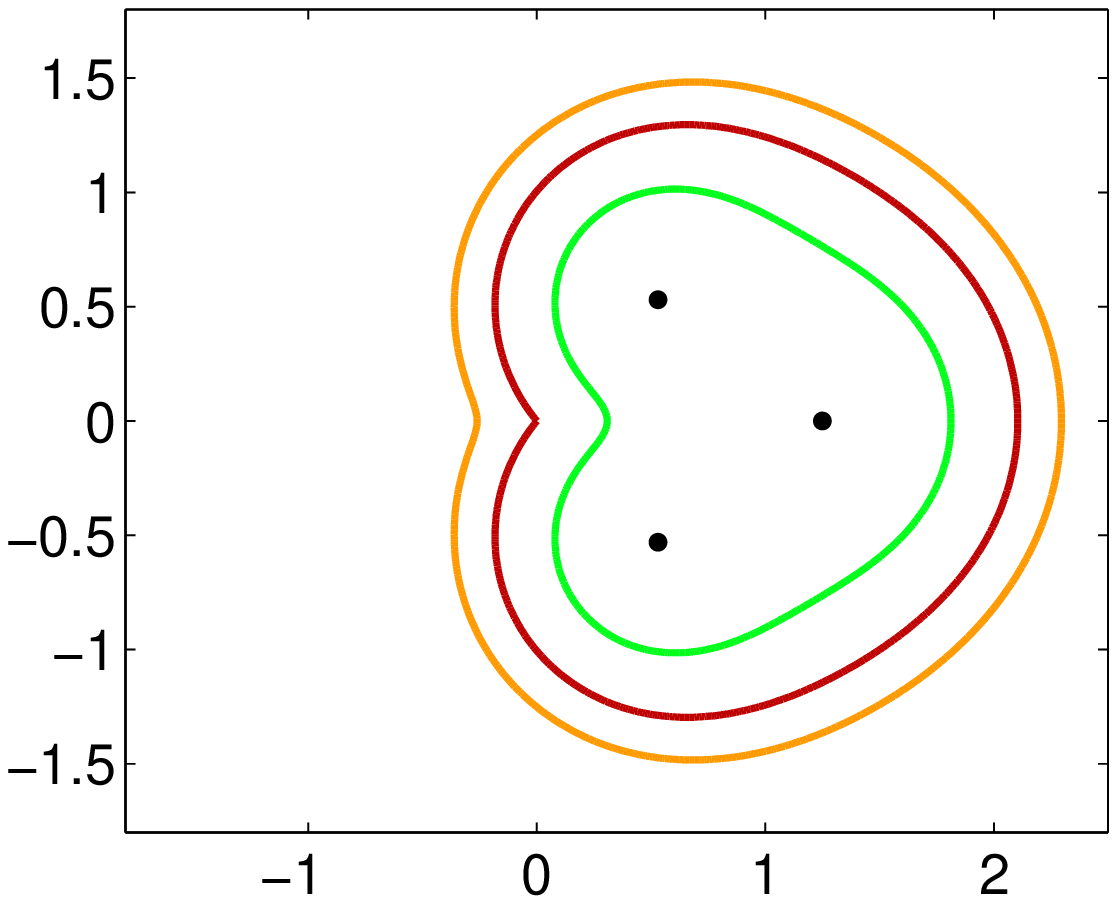, height=4cm} \hspace{-2cm}
\epsfig{file=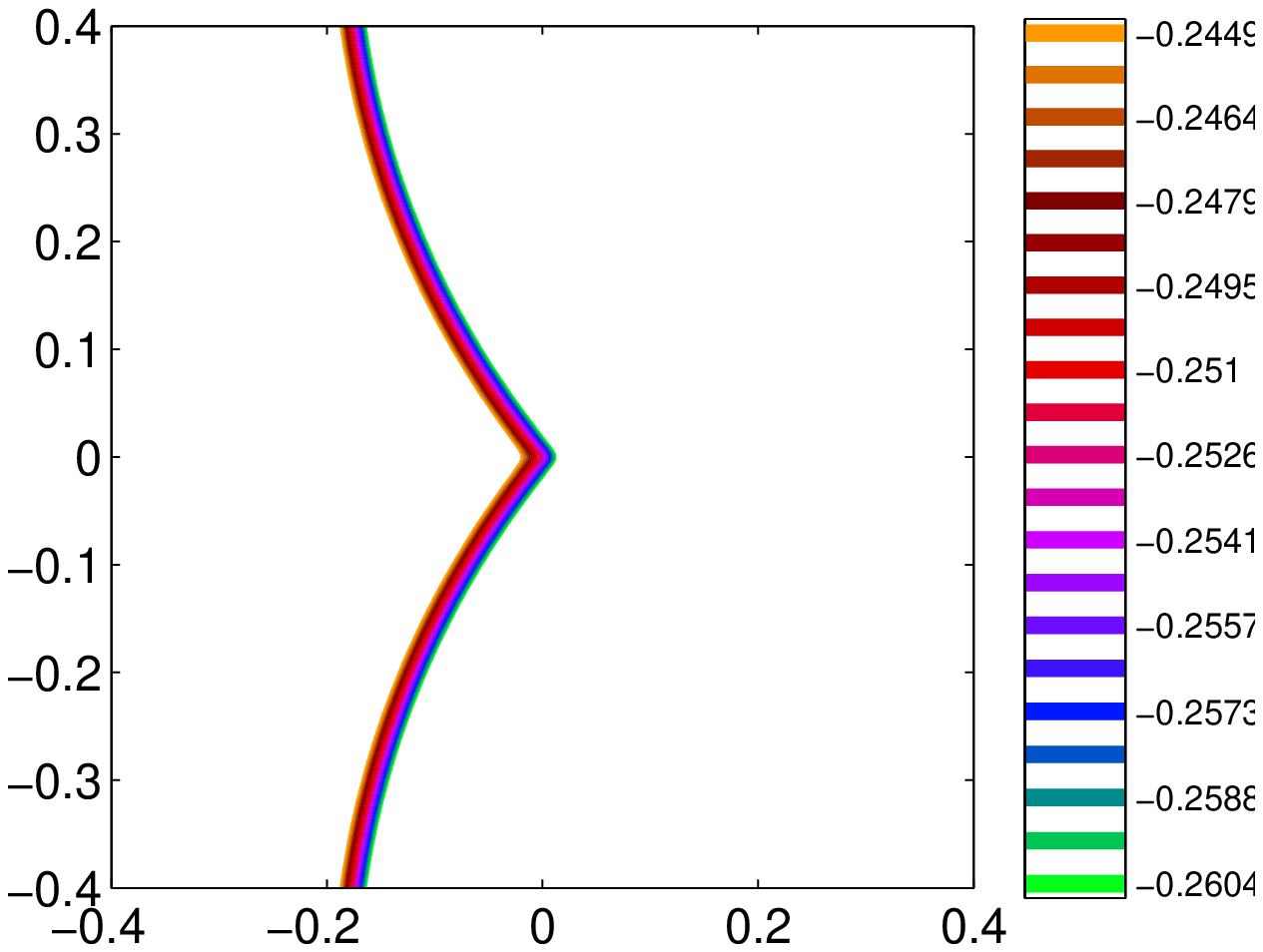, height=4cm}

\epsfig{file=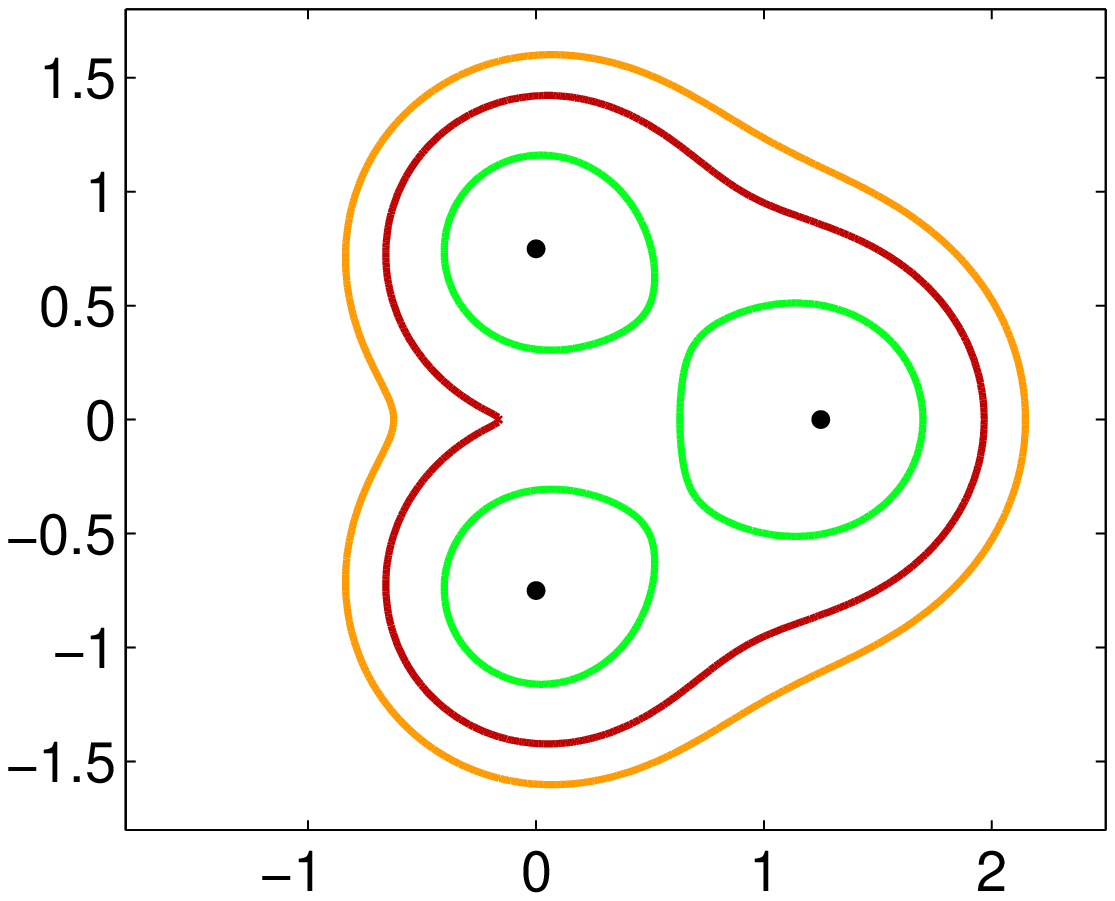, height=4cm} \hspace{-2cm}
\epsfig{file=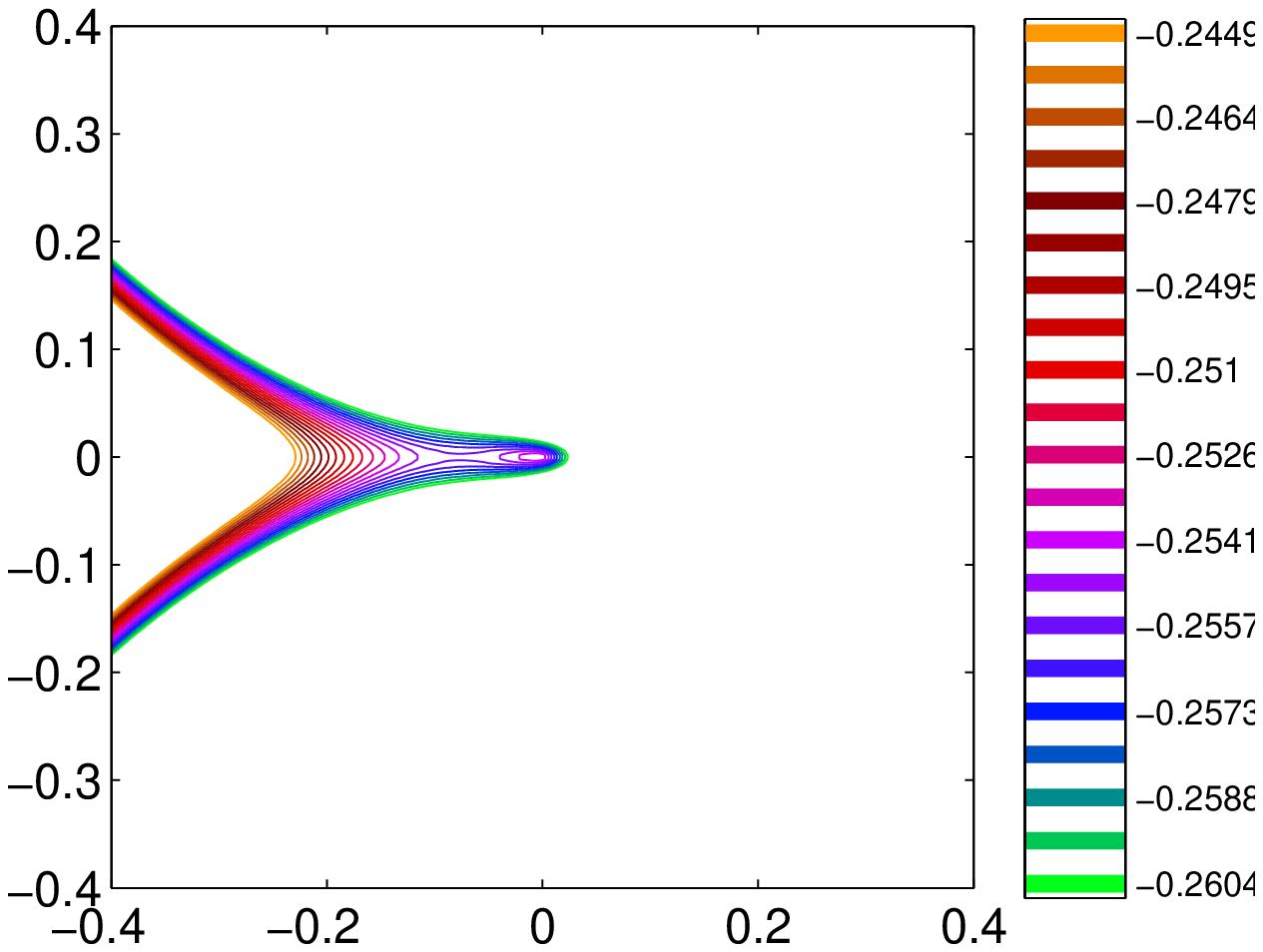, height=4cm}

\epsfig{file=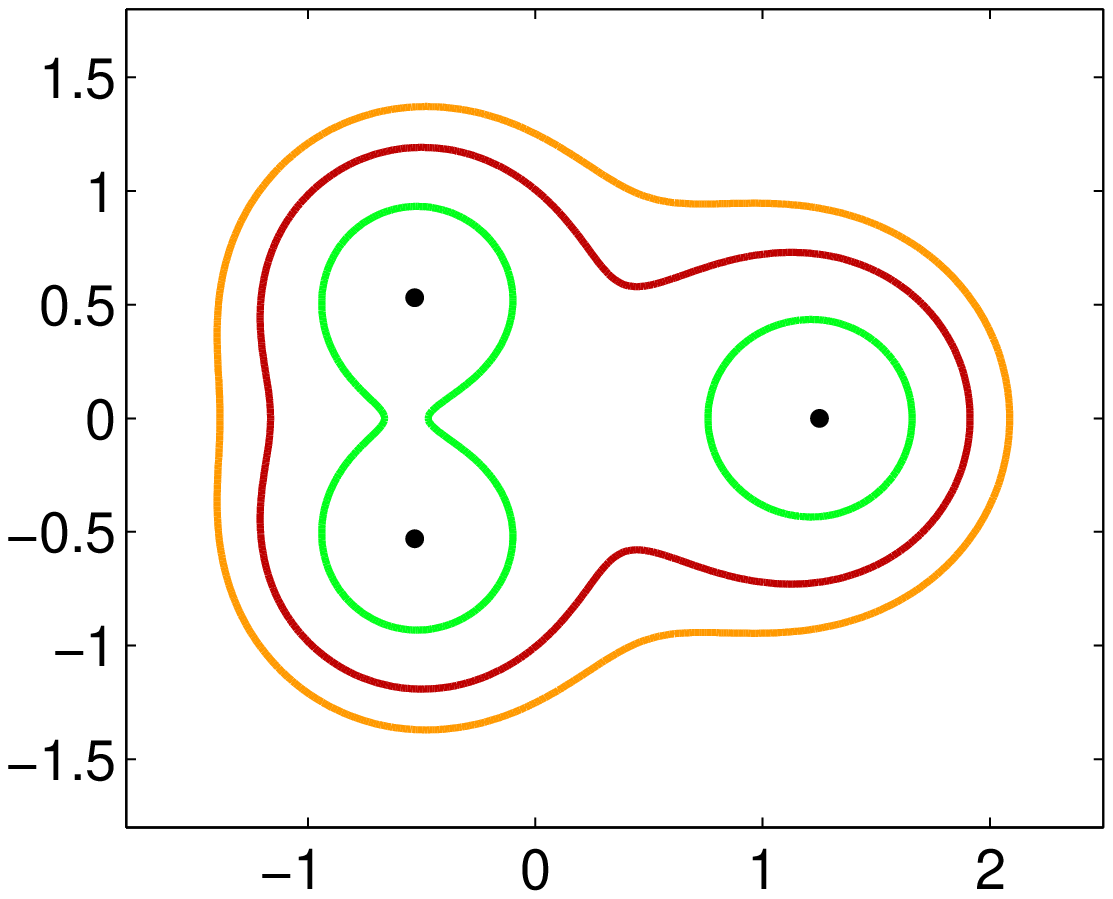, height=4cm} \hspace{-2cm}
\epsfig{file=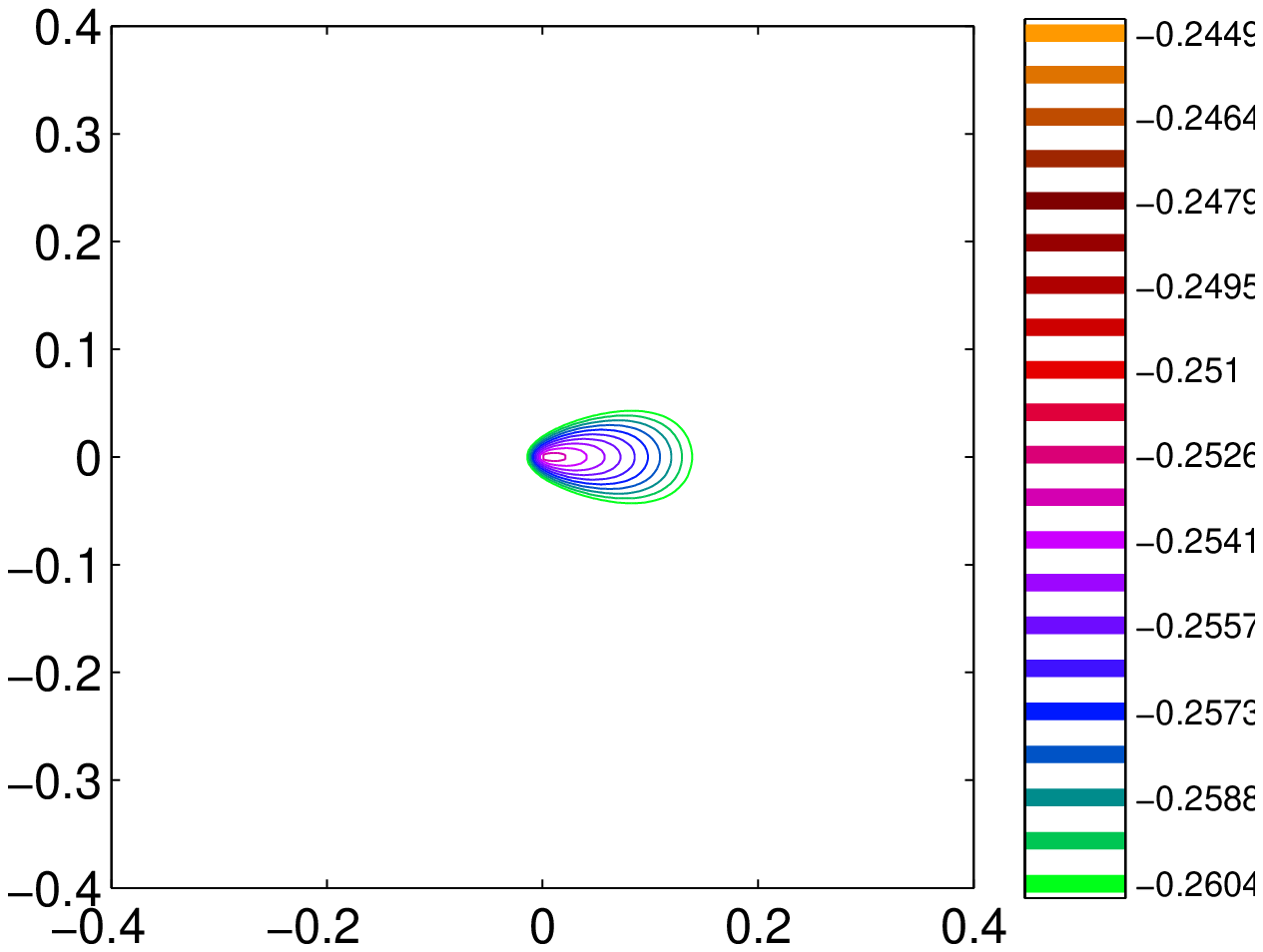, height=4cm}

\epsfig{file=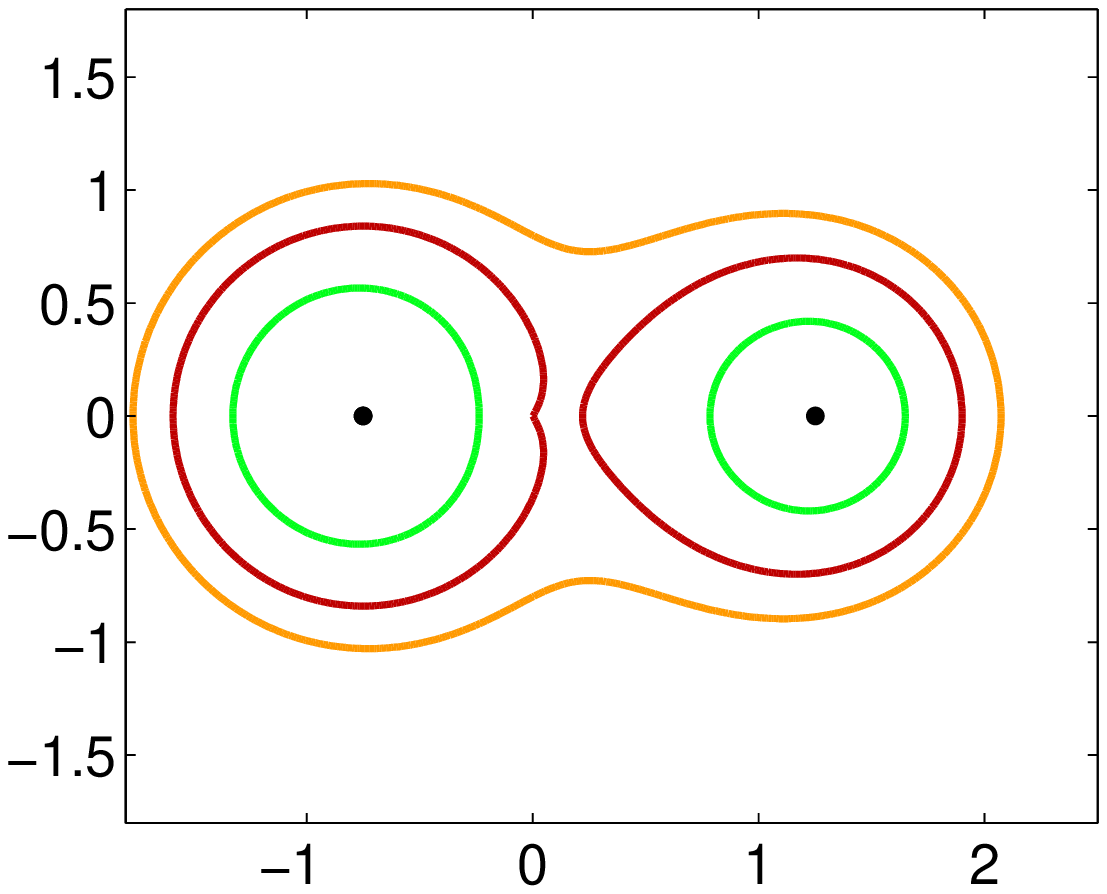, height=4cm} \hspace{-2cm}
\epsfig{file=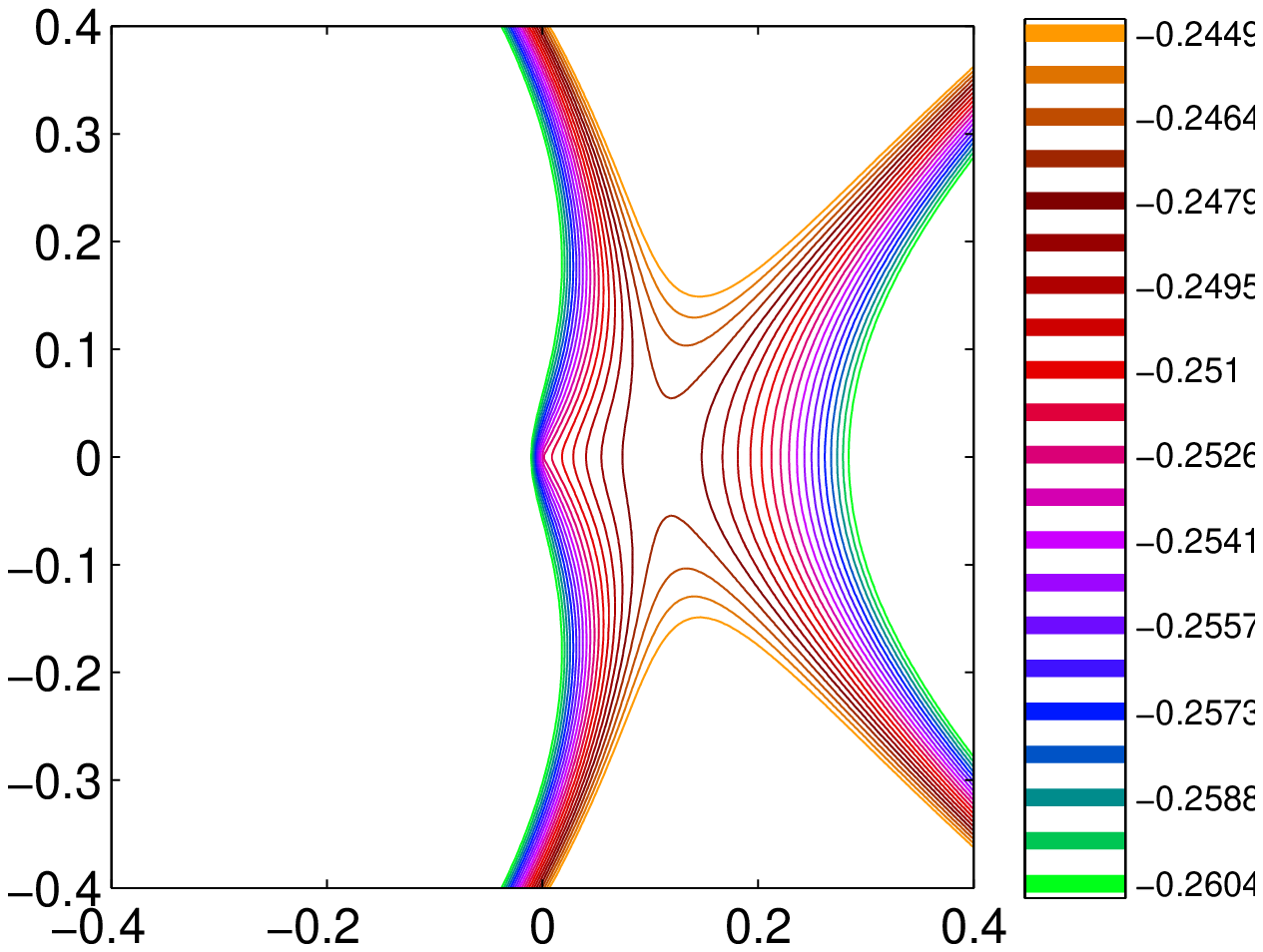, height=4cm}

\caption{Evolution of the essential singularity on the pseudospectral boundary
for the family of matrices given in Example~\ref{ex4}. Here 
$a=\overline{c}=3e^{ik\pi/4}/4$
for $k=0,1,2,3,4$ (top to bottom).  On the left side 
we depict $\partial \spe_\delta A$ for
$\delta=\frac{3}{4\sqrt{5}},\frac{\sqrt{5}}{4},\frac{41}{4\sqrt{5}}$.
On the right side we show details of the pseudospectral 
boundaries near 
the origin for $\delta$ close to the critical value $\sqrt{5}/4$.
\label{f2}}

\end{figure}


\begin{thebibliography}{99}

\bibitem{abo}
{\scshape R. Alam, S. Bora}, On sensitivity of eigenvalues and
eigendecomposition of matrices, \textit{Linear Algebra Appl.},
\textbf{396} (2005), pp. 273-301.

\bibitem{blp}{\scshape L.~Boulton, P.~Lancaster, P.~Psarrakos},
On pseudospectra and their boundaries, To appear in 
\emph{Mathematics of Computation} (2007).

\bibitem{fiedler}{\scshape M.~Fiedler},
\emph{Special Matrices and their Applications in Numerical Analysis}, 
Martinus Nijhoff Publisher, 1986.

\bibitem{glr}{\scshape I.~Gohberg, P.~Lancaster, L.~Rodman},
\emph{Matrix Polynomials}, Academic Press, 1982.

\bibitem{gol}{\scshape G.~Golub, C.~van Loan}, \emph{Matrix Computations},
North Oxford Academic Publishing, 1986.

\bibitem{k}{\scshape T.~Kato},
\emph{Perturbation Theory of Linear Operators}, Springer, 1980.

\bibitem{lit} {\scshape D.E.~Littlewood}, On unitary equivalence,
\textit{J. London Math. Soc.}, \textbf{28} (1953), pp. 314-322.

\bibitem{sha} {\scshape H.~Shapiro},
A survey on canonical forms and invariants for unitary
similarity, \textit{Linear Algebra Appl.},
\textbf{147} (1991), pp. 101-167.

\bibitem{tre} {\scshape L.N.~Trefethen, M.~Embree}, 
\emph{Spectra and Pseudospectra}, Princeton University press, 2005.

\bibitem{wil}{\scshape J.H.~Wilkinson},
\emph{The Algebraic Eigenvalue Problem}, Clarendon Press, 1965.

\end{thebibliography}
\end{document}